\documentclass[a4paper,11pt,miktex,win]{samgu}
\usepackage{amssymb, amsfonts, amsthm, array}
\usepackage{mathdots}

\usepackage{epsf}

\runtitle{Field of Borel Group Invariant of adjoint representatnion of the group GL(n, K)}
\runauthor{К.А.\,Вяткина, А.Н. \,Панов}

\begin{document}

\initcounters
\beforearticle[MATH]

\begin{center}
\language=0
\maintitle[Communicated by Postgraduate of Dept. of Algebra and Geometry Vyatkina K.A.]
{FIELD OF BOREL GROUP INVARIANT OF ADJOINT REPRESENTATION OF THE GROUP GL(n,K)}
\authorright[2014]
\anauthor[Vyatkina Kseniya Anatoljevna
\email{vjatkina.k@gmail.com}, Dept. of Algebra and Geometry,
Samara State University, Samara, 443011, Russia.]{K.A.\,Vyatkina}{}
\end{center}

\begin{abstract}
\language=0
The paper is devoted to invariant theory problems. In particular, to the problem of finding generators of invariant fields in an explicit form. The set of generators is given for invariant field of unitriangular group of adjoint representation of $GL(n,K)$. Moreover, the set of generators of borel group invariant field is constructed and their algebraic independence is proved.       
\end{abstract}

\section{Preliminaries}
The main problem of the theory of invariants of algebraic groups is a description of the structure of the ring and the field of invariants. There are two general theorem. \\


{\bf Theorem 1.1.} (\cite{M})
\label{Th_M}
Let $V$ be a vector space over the field $K$. The field of invariants with respect to any triangularizable subgroup of $GL(V)$ is rational. It means, that there are invariants  $Q_1, \dots Q_n \in K(V)$ such that the field of invariants is the field of rational functions $K(Q_1, \dots, Q_n)$.\\


{\bf Definition 1.} A rational transformation $\alpha$ of $K^n$ is triangular if in a suitable coordinate system all transformations $\alpha(g), g\in G$, have the form of
$$ x_i \rightarrow a_ix_i +f_i(x_{i+1}, \dots, x_n),$$
where $a_i \in K^*, f_i \in k(x_{i+1}, \dots , x_n)$.\\


{\bf Theorem 1.2.} (\cite{V})
\label{Th_V} Let $K$ be an algebraically closed field of characteristic zero. For any subgroup, consisting of rational triangular transformations $ K^n$, a field of invariants is rational.\\


Now introduce some notation. Let $K$ be a field of characteristic zero, $V$ be a vector space over the field $K$, $G = \textmd{GL}(n,K)$ is a group of nonsingular matrices, and $B$ is uppertriangular subgroup and $U$ is upperunitriangular subgroup of $G$.
We consider the adjoint action of the group $G$ in the vector space $V
= \textmd{Mat}(n,K) ~\colon ~Ad_g X=gXg^{-1},~  g\in G,~ X \in V$.
 Let's define representation $Ad_g$ in $K[V]$ (accordingly, $\mathcal{F}=K(V)$) by the formula:
$$ Ad_g f(x) = f( Ad_g^{-1} ( x)).$$
Denote $\mathcal{F}^U$ (accordingly $\mathcal{F}^B$) as subfield of $U$-invariants (accordingly, $B$-invariants) in the field $\mathcal{F}$. From theorems 1.1. 
and 1.2.
the fields $\mathcal{F}^U$ and $\mathcal{F}^B$ are rational.


The main goal of this work is to find generators of the field $B$-invariants for the adjoint representation in explicit form. The main result will be formulated in the theorem 2.2.

\section{The main results.}
\subsection{Generators of $U$-invariants field.}
\label{firstJ_k}
Let $x_{i,j}$ be a system of standard coordinate functions on $V=\textmd{Mat}(n,K)$. Let's construct a system of $n$ minors of the form:
\begin{equation}
\label{J_i}
J_{i}=\left| \begin{array}{ccc}
x_{n-i+1,1} & \cdots &x_{n-i+1,i}\\
\vdots & \cdots & \vdots \\
x_{n, 1} & \cdots & x_{i, i}\\
\end{array} \right| ,  i=\overline{1,n}.
\end{equation}

Extend this set to the set of invariants of $ \frac{n (n +1)}{2}$ elements as follows.
We associate with each  $J_{i}$ the system of $i$ determinants $J_{i,j}$, where
$0\leq j\leq i-1$. The determinant $J_{i,0}$ coincides with the minor $J_i$.
The first $i-j$ rows in the determinant $J_{i,j}$, where $1\leq j\leq i-1$, 
coincide with the last $i-j$ rows in the minor $J_i$; the last $j$ rows in $J_{k,i}$ coincide with the similar rows in the minor $J_i(X^*)$
of the adjoint matrix $X^* = (x_{i,j}^*)$, which is defined as:
$$X\cdot X^*= X^* \cdot X = \det X \cdot E$$

{\bf Exapmle 1}. The case of $ n=2$,\quad

$$X=\left(\begin{array}{cc}
x_{11}&x_{12}\\
x_{21}&x_{22}\end{array}\right),\quad  X^*=\left(\begin{array}{cc}
x^*_{11}&x^*_{12}\\
x^*_{21}&x^*_{22}\end{array}\right),$$\\
$ J_{1,0} = x_{2,1}$,\\
$J_{2,0} = \left|\begin{array}{cc}
x_{11}&x_{12}\\
x_{21}&x_{22}\end{array}\right|$,\quad $J_{2,1} =
\left|\begin{array}{cc}
x_{21}&x_{22}\\
x^*_{21}&x^*_{22}\end{array}\right|$.\\

{\bf Example 2}.The case of $\:n=3,$

$$\:X=\left(\begin{array}{ccc}
x_{11}&x_{12}&x_{13}\\
x_{21}&x_{22}& x_{23}\\
x_{31}&x_{32}& x_{33}\end{array}\right),\:X^*=\left(\begin{array}{ccc}
x^*_{11}&x^*_{12}&x^*_{13}\\
x^*_{21}&x^*_{22}& x^*_{23}\\
x^*_{31}&x^*_{32}& x^*_{33}\end{array}\right),$$\\
$ J_{1,0} = x_{3,1}$,\\
$J_{2,0} = \left|\begin{array}{cc}
x_{21}&x_{22}\\
x_{31}&x_{32}\end{array}\right|$,\quad $J_{2,1} =
\left|\begin{array}{cc}
x_{31}&x_{32}\\
x^*_{31}&x^*_{32}\end{array}\right|$,\\
\\
$J_{3,0} = \left|\begin{array}{ccc}
x_{11}&x_{12}&x_{13}\\
x_{21}&x_{22}& x_{23}\\
x_{31}&x_{32}& x_{33}\end{array}\right|$,\quad $J_{3,1} =
\left|\begin{array}{ccc}
x_{21}&x_{22}&x_{23}\\
x_{31}&x_{32}& x_{33}\\
x^*_{31}&x^*_{32}& x^*_{33}\end{array}\right|$,\quad $J_{3,2} =
\left|\begin{array}{ccc}
x_{31}&x_{32}&x_{33}\\
x^*_{21}&x^*_{22}& x^*_{23}\\
x^*_{31}&x^*_{32}& x^*_{33}\end{array}\right|$.\\
\\

{\bf Theorem 2.1.}{\bf (Vyatkina K.A, Panov A.N. \cite{VP})} The field of $U$-invariants of adjoint representation of the group 
$\mathrm{GL}(n,K)$ is the field of rational functions of $\{J_{i,j}:\,
1\leqslant i\leqslant n,\,\, 0\leqslant j\leq i-1\}$.
\subsection{Generators of $B$-invariants field.}

Let $H$ denote the subgroup of diagonal matrix in $\mathrm{GL}(n,K)$. Its easy to see, that adjoint representation of subgroup $H$  preserves the field $\mathcal{F}^U$, and the field $\mathcal{F}^B$ coincides with the subfield of $H$-invariants in $ \mathcal{F}^U $.\\

{\bf Definition 2.} We call function $f \in k[V]$ cumulant, if $$Ad_h f = \chi_h (f) \cdot f,$$ where $h \in H$ and  $\chi_h(f) \in K^*$. The character $\chi_h(f)$ of subgroup $H$ we call weight of cumulant $f$.

Let $h=\textmd{ diag} (a_1, \dots, a_n) \in H$. Then $Ad_h x_{i,j} = \displaystyle\frac{a_j}{a_i} x_{i,j}$, that is $$\chi_h (x_{i,j}) = \displaystyle\frac{a_j}{a_i}.$$

 Note, that adjugate matrix $X^*$ changing  similary. Indeed, the definition of the adjoint matrix gives $$ X\cdot X^* = X^* \cdot X = \det X \cdot E,$$
$$Ad_h(X) \cdot Ad_h(X^*)= \det (Ad_h(X)) \cdot E. $$
Thence $Ad_h(X^*)=(Ad_h(X))^*$ and hence $Ad_h x_{i,j}^* = \displaystyle\frac{a_j}{a_i} x^*_{i,j}$. Then \begin{equation}\label{eqeq} \chi_h(x^*_{i,j})=\chi_h(x_{i,j})=\displaystyle\frac{a_j}{a_i}.\end{equation}\\

Let's define some notations: \begin{enumerate}
\item \label{y} $y_1= J_{1,0},~
y_i = \displaystyle \frac{J_{i,i-1}}{J_{i-1,0}} ~\textmd{, где}~  2 \leqslant  i \leqslant  n$,
\item \label{Y} $Y_{i,j}= \displaystyle\frac{J_{i,j}\cdot  y_{n-i+j+1} }{J_{i-1,j}\cdot y_{i}}
~\textmd{, где}~  2  \leqslant  i \leqslant  n,~ 0  \leqslant  j \leqslant  i-2 $.
\end{enumerate}
Now let's prove the next lemma.\\

{\bf Lemma 1.} Let  $F=K(z_1, \dots z_n)$ --- field of a rational functions of variables $z_1,\ldots,z_n$ from field $K$. Let diagonal subgroup $H$ act on the field $F$ according to the formula: $hz_i = \chi_h(z_i)z_i$. Then subfield of $H$-invariants $F^H$ generated by the elements of the form $$z_1^{m_1}\dots z_n^{m_n}\textmd{, where } m_i \in \mathbb{Z} \textmd{ and }{\chi_h(z_1)}^{m_1} \dots {\chi_h(z_n)}^{m_n}=1.$$

{\bf Definition 3.} We call the elementary transformation of a system of generators  $z_1, \dots, z_n$ of the field $F$  the transition from the old system to the new system, using the operations like:
$$z_{i} \rightarrow z_{i}\cdot z_{j}^{\pm 1}\textmd{, where } i\not= j.$$

{\bf Definition 4.} We claim that two system of generators are equivalent, if there are elementary transformation, which turn one system to another.\\

Now let proceed to the statement of the main theorem and prove it.\\

{\bf Theorem 2.2.} The field of $B$-invariants of adjoint representation of the group $\textmd{GL}(n, K)$ is the field of rational functions of
\begin{equation}
\label{L_J} \{y_{n}, \textmd{  } Y_{i,j}:~ 2  \leqslant  i \leqslant  n,~ 0  \leqslant  j \leqslant  i-2 \}
\end{equation}

{ \bf Proof.}
As we say above, $\mathcal{F}^U$ is generated by the set of algebraically independent elements $\{ J_{1,0}, \cdots, J_{n, n-1} \}$ (see Theorem 2.1). We refer to the following table as a table of generators:
\begin{equation}
\label{tab_J}
\left[ \begin{array}{ccccc}
 &  & &  & J_{n,0}\\
 &  &  & J_{n-1,0} & J_{n, 1}\\
 &  & \iddots & \vdots & \vdots\\
 & J_{2,0} &\dots & J_{n-1, n-3} & J_{n,n-2} \\
J_{1,0} & J_{2,1} & \dots & J_{n-1, n-2} & J_{n,n-1} \end{array} \right] 
\end{equation}
Each of generators $J_{i,j}$  is semi-invariant. Similarly to the table of generators we construct the table of weights:
\begin{equation}
\label{tab_chi} \left[ \begin{array}{ccccc}
 &  & &  & {\chi_h}(J_{n,0})\\
 &  & & {\chi_h}(J_{n-1,0}) & {\chi_h}(J_{n, 1})\\
 &  & \iddots & \vdots & \vdots\\
 & {\chi_h}(J_{2,0}) &\dots &  {\chi_h}(J_{n-1,n-3}) & {\chi_h}(J_{n,n-2}) \\
{\chi_h}(J_{1,0}) & {\chi_h}(J_{2,1}) & \dots & {\chi_h}(J_{n-1, n-2}) & {\chi_h}(J_{n,n-1}) \end{array} \right] 
\end{equation}

The action of subgroup $H$ on the field $\mathcal{F}^U$ satisfies the conditions of Lemma 1.  Starting with the system $\{ J_{i,j} \}$ we shell construct by elementary transformation the new system of generators, which has the simplest weighting table, using only elementary transformations.

Apply the formula (\ref{eqeq}),  we have
\begin{equation}
\label{chi}
{\chi_h}(J_{i,j})=\frac{a_{1}\dots a_i }{a_{n}\dots a_{n-j+1}\cdot a_{n}\dots a_{n-i+j+1}}\end{equation}\\
In particular, \\
$${\chi_h}(J_{1,0}) = \displaystyle\frac{a_1}{a_n},\dots,~ {\chi_h}(J_{i, 0}) = \displaystyle\frac{a_1 \dots a_{i}}{a_n \dots a_{n-i+1}},\dots,~ {\chi_h}(J_{n,0})=1,$$
$${\chi_h}(J_{2,1})=\displaystyle\frac{a_1 a_2}{a_n a_n},\dots, ~ {\chi_h}(J_{i,i-1})=\displaystyle\frac{a_1 \dots a_i}{a_n \dots a_{n-i+2} a_n },\dots ,~{\chi_h}(J_{n,n-1})=\displaystyle\frac{a_1}{a_n}.$$
After substitution the obtained values to (\ref{tab_chi}), we obtain:
$$\left[ \begin{array}{cccccc}
 &  & &  & & 1\\
 &  &  &  &\displaystyle\frac{a_1 \dots a_{n-1}}{a_n \dots a_2} & {\chi_h}(J_{n, 1})\\
 &  & & \iddots & \vdots & \vdots\\
 &  & \displaystyle\frac{a_1 a_2 a_3}{a_n a_{n-1} a_{n-2}} & \dots & {\chi_h}(J_{n-1, n-4}) & {\chi_h}(J_{n, n-3})\\
 &  \displaystyle\frac{a_1 a_2}{a_n a_{n-1}} & {\chi_h}(J_{3,1}) & \dots & {\chi_h}(J_{n-1, n-3}) & {\chi_h}(J_{n,n-2}) \\
\displaystyle\frac{a_1}{a_n} & \displaystyle\frac{a_1 a_2}{a_n a_n} & \displaystyle\frac{a_1 a_2 a_3}{a_n a_n a_{n-1}} &\dots & \displaystyle\frac{a_1 \dots a_{n-1}}{a_n a_n a_{n-1} \dots a_3} & \displaystyle\frac{a_1}{a_n}                                                                                                                                                                                                                                                                                                                                                                       \end{array} \right] $$
Let's make the follow series of elementary transformations.

 \begin{enumerate}
\item We proceed from system  $\{J_{i,j}\}$ to system $\{J'_{i,j}\}$, where $ J'_{i,j} = \displaystyle\frac{J_{i,j}}{J_{i-1,j}}$ for $3  \leqslant  i \leqslant  n,~ 1  \leqslant  j \leqslant  i-2$; and $J'_{i,j} = J_{i,j}$ for others.  Since $${\chi_h} \left(J'_{i,j} \right)=\frac{a_i}{a_{n-i+j+1}},$$ the table of weights for $\{ J'_{i,j}\}$ takes the form:
\begin{equation} \label{chi_tilde} \left[ \begin{array} {cccccc}
 &  & &  & 1\\
 &  &  &  \displaystyle\frac{a_1 \dots a_{n-1}}{a_n \dots a_{2}} & \displaystyle\frac{a_n}{a_2}\\
 &  & \iddots & \vdots& \vdots\\
 & \displaystyle\frac{a_1 a_2}{a_n a_{n-1}} & \dots & \displaystyle\frac{a_{n-1}}{a_{n-1}}  & \displaystyle\frac{a_n}{a_{n-1}} \\
\displaystyle\frac{a_1}{a_n} & \displaystyle\frac{a_1 a_2}{a_n a_n} & \dots & \displaystyle\frac{a_1 \dots a_{n-1}}{a_n a_n a_{n-1}\dots a_{3}}  & \displaystyle\frac{a_1}{a_n}                                                                                                                                                                                                                                                                                                                                                                       \end{array} \right] 
\end{equation}

\item Each $ J'_{i, i-1}$, where $2 \leqslant i \leqslant n$, is replaced by $ y_i = \displaystyle\frac{J'_{i, i-1}}{ J'_{i-1,0}}= \displaystyle\frac{J_{i, i-1}}{ J_{i-1,0}}$. Notice, that $y_1=J'_{1,0}=J_{1,0}$. The table of generators becomes:
\begin{equation}
\label{tab_J_1}
\left[ \begin{array}{ccccc}
 &  & &  & J'_{n,0}\\
 &  &  &  J'_{n-1,0} &  J'_{n, 1}\\
 &  & \iddots & \vdots & \vdots\\
 & J'_{2,0} &\dots &  J'_{n-1,n-3} &  J'_{n,n-2} \\
y_{1} & y_{2} & \dots & y_{n-1} & y_{n} \end{array} \right] 
\end{equation}

We calculate ${\chi_h}(y_i)$:
$$
{\chi_h}(y_1)=\displaystyle\frac{a_1}{a_n}, ~\dots,~ {\chi_h}(y_i)=\displaystyle\frac{a_{i}}{a_n},~ \dots~ , {\chi_h}(y_{n-1})=\displaystyle\frac{a_{n-1}}{a_n},~ {\chi_h}(y_n)=\displaystyle\frac{a_n}{a_n}=1$$
Substituting these weights in (\ref{chi_tilde}), we obtain:
$$\left[ \begin{array}{ccccc}
 &  & &  & 1\\
 &  &  & \displaystyle\frac{a_1}{a_n} & \displaystyle\frac{a_n}{a_2}\\
 &  &  \iddots & \vdots & \vdots\\
 & \displaystyle\frac{a_1 a_2}{a_n a_{n-1}} & \dots & \displaystyle\frac{a_{n-1}}{a_{n-1}} & \displaystyle\frac{a_n}{a_{n-1}} \\
\displaystyle\frac{a_1}{a_n} & \displaystyle\frac{a_2}{a_n} & \dots & \displaystyle\frac{a_{n-1}}{a_n}  & 1                                                                                                                                                                                                                                                                                                                                                                       \end{array} \right] $$
\item Now we proceed from the system of generators $$\{y_1,\ldots, y_n, J'_{i,j},~ 2\leqslant i\leqslant n,~ 0\leqslant j\leqslant i-2\}$$ to the new system of elements, defined, setting  $ J''_{i,0}=\displaystyle\frac{  J'_{i,0}}{ J'_{i-1,0}}$ и $J''_{i,j} = J'_{i,j}$ for the other generators. The weights can be calculated by the formula:
$${\chi_h} \left( J''_{i,0}\right) = \frac{a_i}{a_{n-i+1}}, ~\textmd{где}~ 2\le i\le n-1.$$
Table of the weights has the form:
$$\left[ \begin{array}{ccccc}
 &  & &  & 1\\
 &  &  &  \displaystyle\frac{a_{n-1}}{a_{2}} & \displaystyle\frac{a_n}{a_2}\\
 &  & \iddots & \vdots & \vdots\\
 & \displaystyle\frac{a_2}{a_{n-1}} & \dots & \displaystyle\frac{a_{n-1}}{a_{n-1}}  & \displaystyle\frac{a_n}{a_{n-1}} \\
\displaystyle\frac{a_1}{a_n} & \displaystyle\frac{a_2}{a_n} & \dots & \displaystyle\frac{a_{n-1}}{a_n }  & 1                                                                                                                                                                                                                                                                                                                                                                       \end{array} \right] $$
\item  Each $J''_{i,j}$, where $ 2\leqslant i\leqslant n,~ 0\leqslant j\leqslant i-2$, we replace by the $$Y_{i,j} = \displaystyle\frac{J''_{i,j} \cdot y_{n-i+j+1}}{ y_i}.$$   
We obtain
$${\chi_h}(Y_{i,j})={\chi_h}\left( \frac{J''_{i,j} \cdot y_{n-i+j+1}}{ y_i}\right)=\displaystyle\frac{ a_i}{a_{n-i+j+1}} \cdot \displaystyle\frac{a_{n-i+j+1}}{a_n} \cdot \displaystyle\frac{a_n}{a_i}= 1,$$
that is, each $Y_{i,j}$ is invariant.

The table of generators and the table of weights has the form:$$ \left[ \begin{array}{ccccc}
&  & &  & Y_{n,0}\\
&  &  &  Y_{n-1 ,0}& Y_{n, 1}\\
&  & \iddots & \vdots & \vdots\\
& Y_{2,0} &  \dots & Y_{n-1,n-3}  & Y_{n, n-2} \\
y_1 & y_2 &  \dots & y_{n-1}  & y_n                                                                                                                                                                                                                                                                                                                                                                        \end{array} \right] 
\quad\quad\quad
\left[ \begin{array}{ccccc}
&  & &  & 1\\
&  &  & 1 & 1\\
&  & \iddots & \vdots & \vdots\\
& 1 & \dots  & 1& 1 \\
\displaystyle\frac{a_1}{a_n} & \displaystyle\frac{a_2}{a_n} & \dots & \displaystyle\frac{a_{n-1}}{a_n } & 1                                                                                                                                                                                                                                                                                                                                                                        \end{array} \right] 
$$
\end{enumerate}
Applying Lemma 1 to the system of generators $$\{y_k, Y_{i,j}:~ 1  \leqslant  k \leqslant  n-1,~ 2  \leqslant  i \leqslant  n,~ 0  \leqslant  j \leqslant  i-2 \},$$ we conclude that the set of algebraically independent elements  $\{y_{n}, \textmd{  } Y_{i,j}:~ 2  \leqslant  i \leqslant  n,~ 0  \leqslant  j \leqslant  i-2 \}$ generates the field  $\mathcal{F}^B$.

\renewcommand{\refname}{References}

\vspace{5mm}
\noindent Paper received 10/{\it IV}/2014.\\
\noindent Paper accepted 10/{\it IV}/2014.
\vspace{5mm}

\end{document}